\documentclass[11pt]{amsart}
\usepackage{geometry}                
\usepackage{graphicx}
\usepackage{amssymb}
\usepackage{epstopdf}
\usepackage{color}
\usepackage{tikz}
\usetikzlibrary{calc,shapes.geometric,decorations.markings,shadows,decorations.fractals}

\newtheorem{thm}{Theorem}[section]
\newtheorem{prop}[thm]{Proposition} 
 
\newtheorem{lemma}[thm]{Lemma}

\theoremstyle{definition}
\newtheorem{example}{Example}
\newtheorem{defn}[thm]{Definition}

\title{Relative Vertex Asphericity}
\author{Jens Harlander and Stephan Rosebrock}

\begin{document}
\begin{abstract} Diagrammatic reducibility DR and its generalization, vertex asphericity VA,  are combinatorial tools developed for detecting asphericity of a 2-complex. Here we present tests for a relative version of VA that apply to pairs of 2-complexes $(L,K)$, where $K$ is a subcomplex of $L$. We show that a relative weight test holds for injective labeled oriented trees, implying that they are VA and hence aspherical. This strengthens a result obtained by the authors in 2017 and simplifies the original proof.
\end{abstract}
\maketitle

\section{Introduction}

A 2-complex $L$ is {\em vertex aspherical} VA, if every combinatorial map from a 2-sphere into $L$ contains a pair of faces with a vertex  in common so that the faces are mapped mirror-wise across this vertex to the same 2-cell of $L$. Vertex asphericity implies topological asphericity.  The closely related concept of {\em diagrammatic reducibility} DR was introduced by Sieradski \cite{S83} in 1983. See also Gersten \cite{Ger87}. 
The weaker notion of VA was first considered by Huck-Rosebrock \cite{HR02}. Recent developments concerning combinatorial asphericity can be found in Barmak-Minian \cite{BaM18} and Blufstein-Minian \cite {BM19}.

Relative vertex asphericity for pairs of 2-complexes $K\subseteq L$ already appeared in a previous article \cite{HaRo17} by the authors, where it was used to establish asphericity of injective LOT-complexes. 
Other and related notions of relative combinatorial asphericity are in the literature. Diagrammatic reducibility for relative presentations was considered by Bogley-Pride  \cite{BP92} in 1992 and has found many applications over the years. See Bogley-Edjvet-Williams \cite{BEW18} for a good overview. Very recently the idea of directed diagrammatic reducibility was introduced and studied by the authors in \cite{HaRo19}.

A powerful method for showing that a 2-complex $L$ is DR is by showing that it satisfies a weight test. This test appeared first in work of Sieradski \cite{S83}, who called it a coloring test, and was later generalized by Gersten \cite{Ger87} who saw it as a combinatorial version of the Gauss-Bonet Theorem. Sieradski allowed only colors (now called weights) $0$ and $1$, whereas any real number could be used for weights in Gersten's version. Wise \cite{Wise} showed that if $L$ satisfies the coloring test then it has non-positive sectional curvature and hence $\pi_1(L)$ is locally indicable. This is not true in the context of Gersten's weight test. In this paper we give a relative weight test for pairs $(L,K)$ that implies relative VA. 

The Whitehead conjecture, which states that subcomplexes of aspherical 2-complexes are aspherical, has motivated the development of combinatorial versions of asphericity, such as DR and VA, and labeled oriented trees have long been known to be an important testing ground for the conjecture. In \cite{HR01} Huck and Rosebrock proved that prime injective labeled oriented trees satisfy Sieradski's coloring test. In \cite{HaRo17} the authors showed that injective labeled oriented trees are aspherical. We strengthen this result here by showing that injective labeled oriented trees satisfy a relative weight test with weights $0$ and $1$ and hence are VA. The proof is a direct generalization of the proof given in the prime case.\\

We thank the referee for several helpful comments and suggestions.

\section{Relative Vertex Asphericity}

A map $f\colon X\to Y$ between complexes is {\em combinatorial} if $f$ maps open cells of $X$ homeomorphically to open cells of $Y$.  Here a 2-complex will always mean a 2-dimensional cell complex with combinatorial attaching maps. A {\em spherical diagram} over a 2-complex $L$ is a combinatorial map $C\to L$, where $C$ is a 2-sphere with a cell structure. If a 2-complex $L$ is non-aspherical, then there exists a spherical diagram which realizes a
nontrivial element of $\pi_2(L)$. In fact, $\pi_2(L)$ is generated by spherical diagrams. So in order to check whether a 2-complex is aspherical or not it is enough to check spherical diagrams. We also study {\it surface diagrams}. These are combinatorial maps $F\to L$, where $F$ is an orientable surface with or without boundary.

The link of a vertex $u$, lk$(L,u)$, is the boundary of a regular neighborhood of $u$ in $L$. So lk$(L,u)$ is a graph whose edges are the corners of 2-cells at $u$. Suppose $L$ is a standard 2-complex with a single vertex $u$ and oriented edge set $X$. Then the vertices of lk$(L,u)=$ lk$(L)$ are $\{ x^+, x^- \mid x\in X\}$, where $x^+$ is a point of the oriented edge $x$ close to the beginning, and $x^-$ is a point close to the ending of that edge.  The {\it positive link} lk$^+(L)$  is the {\it full subgraph} on the vertex set  $\{ x^+ \mid x\in X\}$ and the {\it negative link} lk$^-(L)$ is the full subgraph on the vertex set $\{x^- \mid x\in X\}$. The positive and negative links were called right and left graph, respectively, in \cite{HR01}. Let $f\colon C\to L$ be a spherical diagram and let $v$ be a vertex in $C$. Restricting to the link we obtain a combinatorial map  $f|_{\mbox{lk}(C,v)}\colon \mbox{lk}(C,v)\to\mbox{lk}(L, f(v))$ for every vertex $v\in C$ and we let $z(v)=c_1\ldots c_q$ be the image, which is a closed edge path in lk$(L, f(v))$. 

\begin{defn} Let $\Gamma$ be a graph and $\Gamma_0$ be a subgraph. Let $z=e_1...e_q$ be a cycle (closed edge path) in $\Gamma$. We say
\begin{enumerate}
\item $z$ is {\em homology reduced} if it contains no pair of edges $e_i,e_j$ such that $e_i=\bar e_j$ (the bar indicates opposite orientation and we read $z$ cyclically);
\item $z$ is {\em homology reduced relative to $\Gamma_0$} if any pair of edges $e_i,e_j$ such that $e_i=\bar e_j$ is contained in $\Gamma_0$.
\end{enumerate}
\end{defn}

Let $f\colon C\to L$ be a spherical diagram. A vertex $v\in C$ is called a {\it folding vertex} if $z(v)=c_1\ldots c_q\in \mbox{lk}(L, f(v))$ is not homology reduced. In that case the pair of 2-cells $(d_i, d_j)$ of $C$ containing the preimages of $c_i$ and $c_j$, respectively, satisfying $c_i=\bar c_j$ is called a {\it folding pair}. We call $f$ {\it vertex reduced} if it does not have a folding vertex. A 2-complex $L$ is called {\it vertex aspherical} VA if each spherical diagram over $L$ has a folding vertex. Clearly VA implies asphericity.

\begin{defn}Let $K$ be a subcomplex of the 2-complex $L$.  We say that {\em $L$ is VA relative to $K$} if every spherical diagram $f\colon C\to L$, $f(C)\not\subseteq K$, has a folding vertex with folding pair of 2-cells in $L-K$. 
\end{defn}

We can phrase relative VA also in the following way. $L$ is VA relative to $K$ if in every spherical diagram $f\colon C\to L$, $f(C)\not\subseteq K$, there is a vertex $v\in C$ so that $z(v)=c_1...c_q\subseteq \mbox{lk}(L)$ is not contained in lk$(K)$ and is not homology reduced relative to lk$(K)$. 

\begin{thm}\label{sVAsub} Let $L$ be a 2-complex and $K$ a subcomplex. If $K$ is VA and $L$ is VA relative to $K$ then $L$ is VA.\end{thm}

\noindent Proof. Assume $f\colon C\to L$ is a vertex reduced spherical diagram. Since $L$ is VA relative to $K$ we have that $f(C)\subseteq K$. So $f\colon C\to K$ is a vertex reduced spherical diagram, contradicting the assumption that $K$ is VA. \qed\\

\begin{thm}\label{sVA} If $L$ is VA relative to $K$, then $\pi_2(L)$ is generated, as $\pi_1(L)$-module, by the image of $\pi_2(K)$ under the map induced by inclusion. In particular, if $K$ is aspherical, then so is $L$.
\end{thm}

\noindent Proof. Every vertex reduced spherical diagram $f\colon C\to L$ has its image $f(C)$ in $K$. Thus $f$ represents an element in $\pi_2(K)$. Since $\pi_2(L)$ is generated by vertex reduced spherical diagrams, it follows that $\pi_2(L)$ is generated by the image of $\pi_2(K)$. \qed\\

\section{Tests for relative vertex asphericity}

Let $K\subseteq L$ be 2-complexes. We say a spherical diagram $f\colon C\to L$ is $K$-{\em thin} if for every vertex $v \in C$ there is a 2-cell in $C$ containing $v$ which is mapped to a 2-cell in $L-K$. Thus if $f(\mbox{lk}(C,v))=z(v)=c_1\ldots c_k$, then at least one corner $c_i\in \mbox{lk}(L)-\mbox{lk}(K)$. We can apply a ``reversed subdivision" to $C$ to turn a spherical diagram into a thin one. The idea is to collect material in $C$ that forms an open disc in $C$ which is mapped to $K$ and make it into a single 2-cell in $C$. For this to work we need to attach additional 2-cells to $K$.

\begin{defn}\label{def:thinning} Given a pair $(L,K)$ of 2-complexes, where $K$ is a subcomplex of $L$, we say that $(\hat L, \hat K)$ is a {\em thinning expansion} if the following holds
\begin{enumerate}
\item $\hat K$ is obtained from $K$ by adding 2-cells and $\hat L=L\cup \hat K$;
\item If there exists a vertex reduced spherical diagram $f\colon C\to L$, $f(C)\not\subseteq K$, then there also exists a $\hat K$-thin vertex reduced spherical diagram $f'\colon C'\to \hat L$, $f'(C')\not\subseteq \hat K$.
\end{enumerate}
\end{defn}
Thinning expansions always exist which can be seen in the following example.

\begin{example}\label{ex:min} Let $(L,K)$ be a 2-complex pair and consider a spherical diagram $f\colon C\to L$. 
We can remove open discs from $C$ to obtain a planar diagram $g\colon F\to L-\{\mbox{open 2-cells of $K$}\}$, where $F$ is a connected planar region. Each boundary component $S$ of $F$ maps to $K$ and presents a trivial element in $\pi_1(K)$. For each $S$ we attach a 2-cell to $K$ using $g\colon S\to K$ as attaching map. We do this for all spherical diagrams over $L$ and arrive at a complex $\hat K$. Note that we can attach discs to $F$ and produce a thin spherical diagram $f'\colon C'\to \hat L=L\cup\hat K$. If $C$ does not contain a folding pair $(d_1,d_2)$, $d_i\in L-K$, then neither does $f'\colon C'\to \hat L=L\cup\hat K$.
This construction gives the {\em minimal thinning expansion} of $(L,K)$. Note that $(\pi_1(\hat L), \pi_1(\hat K))=(\pi_1(L), \pi_1(K))$. 
\end{example}

\begin{example} The {\em maximal thinning expansion} is obtained by adding a 2-cell for every closed edge path in $K$. In this case $(\pi_1(\hat L), \pi_1(\hat K))\not =(\pi_1(L), \pi_1(K))$.
\end{example} 

If $K$ is a 2-complex with oriented edges, then an attaching map of a 2-cell $d\in K$ has {\em exponent sum 0}, when traveling along the boundary of $d$ in clockwise direction, one encounters the same number of positive as of negative edges. Here is a setting we will be using for applications.

\begin{example}\label{ex:useful} Suppose $K=K_1\vee\ldots\vee K_n\subseteq L$ and that the attaching maps of 2-cells of $K$ have exponent sum zero. We construct $\hat K_i$ from $K_i$ by attaching 2-cells to every closed edge path in $K_i$ of exponent sum zero. We let $\hat K=\hat K_1\vee\ldots\vee \hat K_n$ and let $\hat L=L\cup \hat K$. Note that $(\hat L, \hat K)$ contains the minimal thinning expansion given in Example \ref{ex:min} and hence is itself a thinning expansion.
\end{example}

Let $\Gamma $ be a graph and $\hat\Gamma=\Gamma_1\cup \ldots \cup \Gamma_n$ be a union of disjoint subgraphs.  We do not assume that $\Gamma$, or the $\Gamma_i$, are connected. We write $\Gamma /\hat\Gamma$ for the graph obtained from $\Gamma$ by collapsing each $\Gamma_i$ to a point.

\begin{defn}\label{dreltree} Let $\Gamma$ be a graph and $\hat \Gamma=\Gamma_1\cup \ldots \cup \Gamma_n$ be a disjoint union as above. We say $\Gamma$ is a {\em forest relative to $\hat\Gamma$} if $\Gamma /\hat\Gamma$ has no cycles. $\Gamma $ is called a {\em tree relative to $\hat\Gamma$} if in addition $\Gamma$ is connected.

\end{defn}

Let $C$ be a cell decomposition of the 2-sphere with oriented 1-cells. A {\it source} in $C$ is a vertex with all its adjacent edges point away from it, and a {\it sink} is a vertex with all its adjacent edges point towards it. A 2-cell $d\in C$ is said to have {\it exponent sum 0} if, when traveling along the boundary of $d$ in clockwise direction, one encounters the same number of positive as of negative edges.
The following theorem is due to Gersten (see \cite{Ger87}):

\begin{thm}\label{sngpg} Let $C$ be a cell decomposition of the 2-sphere with oriented edges, such that all 2-cells have exponent sum 0. Then $C$ contains a sink and a source.
\end{thm}

\noindent Proof. Fix a vertex $v \in C$. If $w$ is a vertex in $C$ define $h(w)$ to be the exponent sum of an edge path in $C$ that connects $v$ to $w$. The height $h(w)$ is well defined because of the exponent sum zero condition of $C$. A vertex of maximal height is a sink, and a vertex of minimal height is a source. \qed\\


A subcomplex $K$ of a 2-complex $L$ is called {\it full}, if for every 2-cell $d\in L$ where all boundary cells are in $K$ we have $d\in K$. If $K=K_1\vee\ldots\vee K_n\subseteq L$ then lk$^+(K)=\mbox{lk}^+(K_1)\cup \ldots \cup \mbox{lk}^+(K_n)$ is a disjoint union, and 
lk$^-(K)=\mbox{lk}^-(K_1)\cup \ldots \cup \mbox{lk}^-(K_n)$ is a disjoint union as well.

\begin{thm}\label{sforest} Let $K=K_1\vee\ldots\vee K_n\subseteq L$. We assume the attaching maps of 2-cells in $L$ have exponent sum $0$, and the $K_i$ are full. If lk$^+(L)$ is a forest relative to lk$^+(K)$ or  lk$^-(L)$ is a forest relative to lk$^-(K)$ then $L$ is VA relative to $K$. Furthermore, the inclusion induced homomorphism $\pi_1(K_i)\to \pi_1(L)$ is injective for every $i=1,\ldots ,n$.
\end{thm}

\noindent Proof. Let us assume that lk$^+(L)$ is a forest relative to lk$^+(K)$. Consider a thinning expansion $(\hat L, \hat K)$ as in Example \ref{ex:useful}. Note that lk$(L)-\mbox{lk}(K)=\mbox{lk}(\hat L)-\mbox{lk}(\hat K)$, and hence lk$^+(\hat L)$ is a forest relative to lk$^+(\hat K)$. We will first show that there is no $\hat K$-thin vertex reduced spherical diagram $f\colon C\to \hat L$, $f(C)\not\subseteq \hat K$. Suppose that there is such a diagram. Since we assumed the exponent sums of attaching maps of 2-cells in $L$ are $0$, the attaching maps of 2-cells in $\hat L$ have this quality as well. Thus $C$ has a source and a sink by Theorem \ref{sngpg}. Let $v\in C$ be a source. Then $z(v)=c_1\ldots c_q\subseteq \mbox{lk}^+(\hat L)$, the image of the link of $v$ in $C$, is homology reduced, and hence contained in lk$^+(\hat K)$ because lk$^+(\hat L)$ is a forest relative to lk$^+(\hat K)$. But this contradicts thinness of $f\colon C\to \hat L$.

Assume $L$ is not VA relative to $K$. Then by Definition \ref{def:thinning} there exists a $\hat K$-thin vertex reduced spherical diagram $f\colon C\to \hat L$, $f(C)$ not contained in $\hat K$. But we just proved that such a spherical diagram does not exist.

Suppose the map $\pi_1(K_i)\to \pi_1(L)$ is not injective for some $i$. Then there exists a vertex reduced Van Kampen diagram $g\colon D\to L$ such that $g(\partial D)$ is a non-trivial element of $\pi_1(K_i)$. Note that $D$ has to contain 2-cells that are not mapped to $K$ because the map $\pi_1(K_i)\to \pi_1(K)=\pi_1(K_1)*\ldots *\pi_1(K_n)$ is injective. The boundary of $D$ has exponent sum zero because the 2-cells in $L$ are attached by maps of exponent sum $0$. We can cap off $D$ with a disc $d_1$ and obtain a spherical diagram $f\colon C\to \hat L$, $f(C)\not\subseteq \hat K$. Note that this spherical diagram is vertex reduced.  If it were not, then there would have to be a folding vertex $v$ on the boundary of $D$ with folding pair $(d_1,d_2)$, where $d_2$ is a 2-cell in $D$. But that would mean that $f(d_1)=f(d_2)$ is a 2-cell in $L$. Since we assumed $K$ is full, this would imply that $f(d_1)$ is a 2-cell in $K_i$, which contradicts the fact that $g(\partial D)$ is a non-trivial element of $\pi_1(K_i)$. By Definition \ref{def:thinning} there exists a $\hat K$-thin vertex reduced spherical diagram $f'\colon C'\to \hat L$, $f(C)\not\subseteq\hat K$. But we know already from the beginning of this proof that no such spherical diagram exists. \qed\\

Let $L$ be a 2-complex. We assign weights (or angles) $\omega(c)\in \mathbb R$ to the corners of the 2-cells and obtain an {\em angled} 2-complex. If $L=S$ is a closed orientable surface we define the curvature of a 2-cell $d\in S$ to be 
$\kappa(d)=\sum_{i=1}^q \omega(c_i)-(q-2)$, where $c_1,\ldots ,c_q$ are the corners in $d$. The curvature at a vertex is defined to be 
$\kappa(v)=2- \sum \omega(c_i)$
where the $c_i$ are the corners at the vertex $v$. The combinatorial Gauss-Bonet theorem says 
\[ \kappa(S)=\sum_{v\in S} \kappa(v) + \sum_{d\in S} \kappa(d)= 2\chi(S).\]
Note that if $g\colon S\to L$ is a surface diagram and $L$ is an angled 2-complex, we can pull back the weights and give $S$ an induced angle structure. The idea behind a weight test is to give conditions on the link of an angled 2-complex that imply $\kappa(S)\le 0$ for every vertex reduced surface diagram $g\colon S\to L$. This in turn implies that there can not exist vertex reduced spherical diagrams $f\colon C\to L$, and hence $L$ is aspherical. 

A cycle in a graph is called {\it reduced} if no oriented edge in the cycle sequence is followed immediately by its inverse. Here is the weight test as defined by Gersten \cite{Ger87} which implies the asphericity of $L$:

\begin{defn}\label{def:wteststand} Let $L$ be a 2-complex. $L$ satisfies the {\it weight test}, if there are weights $\omega(c)\in \mathbb R$ assigned to the corners of the 2-cells of $L$ such that
\begin{enumerate}
\item $\sum_i \omega(c_i)\le q-2$ if $c_1,\ldots ,c_q$ are the corners of a 2-cell of $L$ and
\item if $z=e_1\ldots e_n$ is a reduced cycle in $lk(L)$ then $\sum_{i=1}^n\omega(e_i)\ge 2$.
\end{enumerate}
\end{defn}

We next define a relative weight test. It is coarse but will be sufficient for the applications we have in mind. Assume $K=K_1\vee\ldots\vee K_n\subseteq L$. We assume $L$ contains a single vertex $v$. We define lk$(L,K)$
in the following way: If $y_1,\ldots ,y_l$ are the edges of $K_i$ then we denote by $\Delta(K_i)$ the full graph on the vertices $y_j^{\pm 1}$ of lk$(K_i)$ together with an edge attached at each $y_j^+$ (a loop at that vertex) and at each $y_j^-$. Every pair of vertices in $\Delta(K_i)$ is connected by an edge, and at every vertex we have a loop. For each $i$ we remove lk$(K_i)$ from lk$(L)$ and insert $\Delta(K_i)$ instead. The resulting graph is lk$(L,K)$. Note that if $\hat K$ is any 2-complex obtained from $K$ by attaching 2-cells and $f\colon C\to \hat L=L\cup\hat K$ is a spherical diagram, then $f(\mbox{lk}(C,v))$ yields a cycle in lk$(L,K)$. Let $\Delta(K)=\Delta(K_1)\cup \ldots \cup \Delta(K_n)$, a disjoint union. Further note that lk$^+(L,K)$ is a forest relative to $\Delta^+(K)$ if and only if lk$^+(L)$ is a forest relative to lk$^+(K)$. And similarly for $-$ in place of $+$.

\begin{defn}\label{def:wtest} Assume $K=K_1\vee\ldots\vee K_n\subseteq L$ and we are in the setting of Example \ref{ex:useful}, that is the attaching maps of 2-cells of $K$ have exponent sum zero. We say $L$ {\em satisfies the weight test} relative to $K$ if there is a weight function $\omega$ on the set of edges of lk$(L,K)$, such that
\begin{enumerate}
\item $\sum_i \omega(c_i)\le q-2$ if $c_1,\ldots ,c_q$ are the corners of a 2-cell of $L$ not contained in $K$,
\item if $z=e_1\ldots e_n$ is a homology reduced cycle in lk$(L, K)$ containing at least one corner from  lk$(L,K)-\Delta(K)$, then $\sum_{i=1}^n\omega(e_i)\ge 2$ and
\item $\omega(c)=0$ if $c$ is an edge of $\Delta^{+}(K_i)$ or $\Delta^{-}(K_i)$,\\
  $\omega(c)=1$ if $c$ connects a vertex of $\Delta^{+}(K_i)$ with one of $\Delta^{-}(K_i)$, $i=1,...,n$.
\end{enumerate}
\end{defn}

\begin{lemma}\label{lemma:gruscht} Assume we are in the setting of Definition \ref{def:wtest}. If $z$ is a cycle in $lk(L, K)$ containing at least one corner from  lk$(L,K)-\Delta(K)$ and $z$ is homology reduced relative to $\Delta(K)$, then $\omega(z)\ge 2$.
\end{lemma}

\noindent Proof. Let $z=c_1...c_q$ be a cycle as in the statement of the lemma. If $z$ is homology reduced then $\omega(z)\ge 2$ since we assume Condition (2) holds. If $z$ is not homology reduced then there exists a pair $c_k, c_l$ satisfying $c_l=\bar c_k$. Since we assume $z$ is homology reduced relative to $\Delta(K)$,  $c_k, c_l\in \mbox{lk}(\Delta(K_i))$ for some $i$. Let $z_1=c_1...c_{k-1}c_{l+1}...c_q$ and $z_2=c_{k+1}...c_{l-1}$. Both $z_1, z_2$  are cycles and are homology reduced relative to $\Delta(K)$. Assume first that both $z_1,z_2$ contain at least one corner from  lk$(L,K)-\Delta(K)$. Then by induction of the cycle length we have $\omega(z_i)\ge 2$. Since $c_k, c_l\in \Delta(K)$ both carry weights $\ge 0$ and we have
$$\omega(z)=\omega(z_1)+\omega(z_2)+\omega(c_k)+\omega(c_l)\ge \omega(z_1)+\omega(z_2)\ge 2+2=4.$$
For the remaining case we assume $z_1$ contains a corner from  lk$(L,K)-\Delta(K)$ but $z_2$ does not. But in that case $z_2\subseteq \Delta(K_i)$. Then $\omega(z_1)\ge 2$ and $\omega(z_2)\ge 0$. We have 
\[\omega(z)=\omega(z_1)+\omega(z_2)+\omega(c_k)+\omega(c_l)\ge \omega(z_1)+\omega(z_2)\ge 2.\]
\qed

\begin{thm}\label{thm:WTD} Assume $K=K_1\vee\ldots\vee K_n\subseteq L$, each $K_i$ is full and we are in the setting of Example \ref{ex:useful}, that is the attaching maps of 2-cells of $K$ have exponent sum zero. Assume further that $L$ satisfies the weight test relative to $K$. Then $L$ is VA relative to $K$. If in addition the attaching maps of the 2-cells of $L$ have exponent sum zero, then all the inclusion induced homomorphisms $\pi_1(K_i)\to \pi_1(L)$ are injective.
\end{thm}

\noindent Proof. Let $(\hat L, \hat K)$ be the thinning expansion constructed in Example \ref{ex:useful}. We first make $\hat L$ into an angled 2-complex. Since $\hat L-\hat K=L-K$, we have already weights on the corners of 2-cells in $\hat L-\hat K$. If $\hat d$ is a 2-cell of $\hat K$ we assign to corners in lk$^+(\hat K)$ and in  lk$^-(\hat K)$ weight $0$, and weight $1$ to all the other corners in $\hat d$. 
Suppose $L$ is not $VA$ relative to $K$. Then there exists a vertex reduced spherical diagram $f\colon C\to L$ that is not already a diagram over $K$. By Definition \ref{def:thinning} there also exists a $\hat K$-thin vertex reduced spherical diagram $f'\colon C'\to \hat L$ that is not already a diagram over $\hat K$. We pull back the weights of $\hat L$ and thus turn $C'$ into an angled 2-complex. Condition (1) in the weight test implies that the curvature of a 2-cell not mapped to $\hat K$ is $\le 0$. 
If $d\in C'$ is a 2-cell which is mapped to $\hat K$ it has exponent sum 0, so there are at least 2 corners with weight 0 (the other corners of $d$ have weight 1). So the curvature of $d$ will also be $\le 0$.

Since $f'\colon C'\to \hat L$ is $\hat K$-thin and vertex reduced, for every $v\in C'$ the image $f'(\mbox{lk}(C',v))$ yields a cycle $z\in $ lk$(\hat L,\hat K)$ that is homology reduced relative to $\Delta(K)$ and contains a corner from a 2-cell in $L-K$. Thus by Lemma \ref{lemma:gruscht} $\omega(z)\ge 2$ which implies that the curvature at $v$ in $C'$ is $\le 0$.
So the curvature of $C'$ is $\le 0$. This is a contradiction because $C'$ is a 2-sphere and so the curvature is $2$.

Injectivity of the homomorphisms $\pi_1(K_i)\to \pi_1(L)$ follows by the arguments already provided in the proof of Theorem \ref{sforest}.
\qed\\

A {\em $(\epsilon ,\delta )$-corner} of lk$(L,K)$ for $\epsilon ,\delta \in\{ +,-\}$ is a corner between vertices $x_i^\epsilon $ and $x_j^\delta $ for some $i,j$.

\begin{thm}\label{thm:forest}
Let $K=K_1\vee\ldots\vee K_n\subseteq L$. Assume
\begin{enumerate}
\item the attaching maps of 2-cells in $L$ have exponent sum zero;
\item lk$^+(L)$ is a forest relative to lk$^+(K)=\mbox{lk}^+(K_1)\cup \dots \cup \mbox{lk}^+(K_n)$ and  lk$^-(L)$ is a forest relative to lk$^-(K)=\mbox{lk}^-(K_1)\cup \dots \cup \mbox{lk}^-(K_n)$.\end{enumerate}
With the assignment $\omega\colon \mbox{Edges of lk$(L,K)$}\to \{ 0,1\}$
\begin{align*}
\omega(c) = \left\{ \begin{array}{ll} 
                0 & \hspace{5mm} \mbox{if $c$ is a $(++)$-corner or a $(--)$-corner} \\
                1 & \hspace{5mm}  \mbox{if $c$ is a $(+-)$-corner} \\
                           \end{array} \right.
\end{align*}
$L$ satisfies the weight test relative to $K$.
\end{thm}

\noindent Proof. Note that lk$^+(L)$ is a forest relative to lk$^+(K)$ if and only if lk$^+(L,K)$ is a forest relative to $\Delta^+(K)$ (the same for $-$ in place of $+$). 
Let $(\hat L, \hat K)$ be a thinning expansion as in Example \ref{ex:useful}. Since we assumed that the attaching maps for the 2-cells of $L$ have exponent sum zero, the same is true for the 2-cells of $\hat L$. Thus if $c_1,\ldots ,c_q$ are the corners in a 2-cell of $\hat L$, then there is at least one $(++)$ and one $(--)$-corner among them. So $\sum \omega(c_i)\le q-2$ and the first condition of the weight test holds.

Let $z$ be a homology reduced cycle in lk$(L, K)$ containing at least one corner from a 2-cell of $L-K$. 
If $z$ contains one $(+-)$-corner it has to contain at least two $(+-)$-corners and then $\omega (z)\ge 2$. So assume $z$ contains only $(++)$-corners (or only $(--)$-corners). But since  lk$^+(L,K)$ is a forest relative to $\Delta^+(K)$ and $z$ contains a corner of $L-K$, $z$ is homology reducible. This is a contradiction.
\qed\\

\section{Applications to Labelled Oriented Trees}

A standard reference for labeled oriented graphs, LOG's for short, is \cite{Ro18}. Here are the basic definitions. A LOG is an oriented finite graph $\Gamma$ on vertices $\textbf{x}$ and edges ${\bf e}$, where each oriented edge is labeled by a vertex. Associated with a LOG $\Gamma $ is the {\em LOG-complex} $K(\Gamma)$, a 2-complex with a single vertex, edges in correspondence with the vertices of $\Gamma$ and 2-cells in correspondence with the edges of $\Gamma$. The attaching map of a 2-cell $d_e$ is the word  $xz(zy)^{-1}$, where $e$ is an edge of $\Gamma$ starting at $x$, ending at $y$, and labeled with $z$. 

A labelled oriented graph is called {\em compressed} if no edge is labelled with one of its vertices. A LOG $\Gamma $ is called {\em boundary reducible} if there is a boundary vertex $x\in\Gamma$ which does not occur as an edge label and {\em boundary reduced} otherwise.
A LOG is {\em injective} if each vertex occurs as an edge label at most once. An injective LOG is called {\em reduced} if it is compressed and boundary reduced. A {\em labeled oriented tree}, LOT, is a labeled oriented graph where the underlying graph is a tree. 

It can be shown that an injective LOT can be transformed into a reduced injective LOT without altering the homotopy-type of the LOT-complex. If $\Gamma$ is a LOT and $\Gamma_1$ is a sub-tree of $\Gamma$ with at least one edge, such that each edge label of $\Gamma_1$ is a vertex of $\Gamma_1$, then we call $\Gamma_1$ a {\em sub-LOT} of $\Gamma$. A sub-LOT $\Gamma_1$ of $\Gamma$ is {\it proper}, if $\Gamma_1\ne\Gamma$. A sub-LOT $\Gamma_1$ of $\Gamma$ is {\em maximal} if it is not contained in a bigger proper sub-LOT of $\Gamma$. A LOG is called {\em prime} if it does not contain proper sub-LOTs.

The LOT consisting of a single edge labeled by one of its vertices is prime but not boundary reduced. All other prime LOTs are boundary reduced, since a boundary reducible LOT containing more than one edge contains a proper sub-LOT. It follows that compressed prime LOTs containing at least two edges are boundary reduced and hence reduced.\\

In \cite{HaRo17} the authors have shown that injective LOTs are aspherical. In this section we prove

\begin{thm}\label{thm:general} Let $\Gamma$ be a compressed injective LOT. Then $K(\Gamma)$ is VA.
\end{thm}

 This generalizes a result obtained by Huck and Rosebrock.

\begin{thm}\label{thm:prime}(Huck-Rosebrock \cite{HR01}) Let $\Gamma$ be a reduced injective prime LOT. Then $K(\Gamma)$ satisfies the weight test with weights from $\{ 0, 1\}$. In particular $K(\Gamma)$ is DR (and therefore also VA).
\end{thm}

\begin{figure}[ht]
\centering

\begin{tikzpicture}

\draw(0,0)--(6,0);
\node [below] at (0,0) {$g$}; \fill (0,0) circle (2pt);
\node [below] at (1,0) {$a$}; \fill (1,0) circle (2pt);
\node [below] at (2,0) {$b$}; \fill (2,0) circle (2pt);
\node [below] at (3,0) {$c$}; \fill (3,0) circle (2pt);
\node [below] at (4,0) {$d$}; \fill (4,0) circle (2pt);
\node [below] at (5,0) {$e$}; \fill (5,0) circle (2pt);
\node [below] at (6,0) {$f$}; \fill (6,0) circle (2pt);
\node [above] at (0.5,0) {$f$};
\node [above] at (1.5,0) {$d$};
\node [above] at (2.5,0) {$e$};
\node [above] at (3.5,0) {$b$};
\node [above] at (4.5,0) {$c$};
\node [above] at (5.5,0) {$g$};
\end{tikzpicture}
\caption{\label{fig:nonWT} A reduced injective non-prime LOT which does not satisfy the weight test (with any orientation of its edges). See Huck-Rosebrock  \cite{HR01}.}
\end{figure}
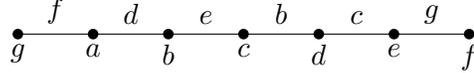

A reduced non-prime injective LOT may not satisfy the weight test. An example is shown in Figure \ref{fig:nonWT}. 

If $\Gamma_1$ is a proper sub-LOT of an injective  LOT $\Gamma $, we can {\it collapse} $\Gamma_1$ in $\Gamma$ and obtain a quotient LOT $\bar\Gamma_1$ in the following specific way:
There exists a unique vertex $x$ in $\Gamma_1$ that does not occur as edge label in $\Gamma_1$, since $\Gamma_1$ is injective. $x$ may or may not occur as an edge label in $\Gamma$. 
Collapse $\Gamma_1$ in $\Gamma$ to $x$ to obtain a quotient tree $\bar\Gamma_1$ of $\Gamma$ with edge set the edges of $\Gamma-\Gamma_1$ and vertex set the vertices that come with the edges in $\Gamma-\Gamma_1$ together with $x$.  If $y$ is the label of an edge $e$ in $\Gamma-\Gamma_1$, then that edge carries the same label $y$ in $\bar \Gamma$. We denote the quotient map by $\Gamma\to \bar \Gamma_1$ and define $\bar z$ to be the image of a vertex $z$ under that map.

Let $\Gamma$ be a reduced injective LOT that is not prime. Choose a maximal proper sub-LOT $\Gamma_1$ and collapse it to obtain the quotient $\bar \Gamma_1$ of $\Gamma$. This quotient is injective. If it is not prime choose a maximal proper sub-LOT $\bar \Gamma_2$ in $\bar \Gamma_1$ and collapse it to obtain $\bar \Gamma_{12}$. We continue this process until we arrive at an injective prime quotient $\bar \Gamma_{12...n}$:

$$\begin{matrix} 
      \Gamma & \to    & \bar \Gamma_1    & \to     &  \bar \Gamma_{12}       & \ldots &\bar\Gamma_{12...(n-1)} &\to & \bar\Gamma_{12...n}\\
      \cup &                & \cup                      &          &   \cup                             &           &   \cup                                                                           \\  
      \Gamma_1 &     &    \bar \Gamma_2 &             & \bar \Gamma_3           &          &\bar\Gamma_n\\
   \end{matrix}$$
Note that $\bar\Gamma_2$ does not contain the vertex $x$ that $\Gamma_1$ is collapsed to, because that would contradict maximality of $\Gamma_1$. Thus the sub-LOT $\Gamma_2$ that maps to $\bar \Gamma_2$ under the first collapse is disjoint from $\Gamma_1$. This procedure gives a disjoint set of sub-LOTs $\Gamma_1,\ldots ,\Gamma_n$ of $\Gamma $ so that 
$$\Gamma_i \to \bar \Gamma_i \ \mbox{is a maximal sub-LOT of }\ \bar \Gamma_{12...(i-1)}.$$
We call a set of sub-LOTs  $\Gamma_1,\ldots ,\Gamma_n$ of $\Gamma $ obtained in this fashion {\it complete}, if $\bar\Gamma_{12\ldots n}$ is a compressed injective prime LOT that is not just a vertex. Such a set may or may not exist. It could happen that all injective prime quotients $\bar\Gamma_{12..n}$ that we obtain are not compressed.\\

\noindent A labelled oriented graph $\Gamma'$ is a {\em reorientation} of a labelled oriented graph $\Gamma$ if $\Gamma'$ is obtained from $\Gamma$ by changing the orientation of some edges.
The following lemma was used in the proof of Theorem \ref{thm:prime}. It is not explicitly stated as a lemma in \cite{HR01}, but the first three lines of Section 3, page 288, of \cite{HR01} state the result, and the proof is given on pages 289 and 290.

\begin{lemma}\label{lem:HRreorient}(Huck-Rosebrock \cite{HR01}) If ${\Gamma}$ is a reduced injective prime LOT, then there is a reorientation ${\Gamma'}$ of ${\Gamma}$ such that $lk^+(K(\Gamma'))$ and $lk^-(K(\Gamma'))$ are trees. In particular $K(\Gamma')$ satisfies the weight test by assigning weight $0$ to all corners in  $lk^+(K(\Gamma'))$ and $lk^-(K(\Gamma'))$, and weight $1$ to all other corners.
\end{lemma}

We can adapt this to our more general setting:

\begin{lemma}\label{lem:reor} If ${\Gamma}$ is a reduced injective LOT with a complete set of sub-LOTs $\Gamma_1,\ldots ,\Gamma_n$, then there is a
reorientation ${\Gamma'}$ of ${\Gamma}$ such that $lk^+(K(\Gamma'))$ and $lk^-(K(\Gamma'))$ are trees relative to $lk^+(K(\Gamma_1'\cup\ldots\cup\Gamma_n'))$ and $lk^-(K(\Gamma_1'\cup\ldots\cup\Gamma_n'))$, respectively. In particular $K(\Gamma')$ satisfies the relative weight test by assigning weight $0$ to all corners in $lk^+(K(\Gamma'))$ and in $lk^-(K(\Gamma'))$, and weight $1$ to all other corners.
\end{lemma}

\noindent Proof. Collapsing each $\Gamma_i$ in $\Gamma$ results in an injective compressed prime LOT $\bar \Gamma =\bar\Gamma_{1...n}$. By Lemma \ref{lem:HRreorient} we can reorient $\bar \Gamma$ to $\bar \Gamma'$ so that both $lk^+(K(\bar \Gamma'))$ and $lk^-(K(\bar \Gamma'))$ are trees. We pull back the edge-orientations of $\bar \Gamma'$ to edge-orientations of $\Gamma $ to achieve a reorientation $\Gamma'$ of $\Gamma$. Note that this reorientation does not affect the $\Gamma_i$ (so $\Gamma_i'=\Gamma_i$). Since both $lk^+(K(\bar \Gamma'))$ and $lk^-(K(\bar\Gamma'))$ are trees, $lk^+(K(\Gamma'))$ is a tree relative to $lk^+(K(\Gamma_1'\cup\ldots\Gamma_n'))$ and $lk^-(K(\Gamma'))$ is a tree relative to $lk^-(K(\Gamma_1'\cup\ldots\Gamma_n'))$. Then Theorem \ref{thm:forest} implies that $K(\Gamma')$ satisfies the weight test.\qed

\smallskip The next lemma was also used in the proof of Theorem \ref{thm:prime}:

\begin{lemma}\label{lem:HRweights}(Huck-Rosebrock \cite{HR02}) Let $\Gamma$ be a reduced injective LOT that satisfies the weight test with weights $0$ and $1$. Then any reorientation of $\Gamma $ satisfies the weight test with weights $0$ and $1$.
\end{lemma}

The analogous result for the general situation is 

\begin{lemma}\label{lem:weights} Let $\Gamma$ be a reduced injective LOT with a complete set of sub-LOTs $\Gamma_1,\ldots ,\Gamma_n$. If $K(\Gamma)$  satisfies the weight test relative to $K(\Gamma_1\cup\ldots\cup\Gamma_n)$ with weights $0$ and $1$, then for any reorientation $\Gamma'$ the 2-complex $K(\Gamma')$  satisfies the weight test relative to $K(\Gamma'_1\cup\ldots\cup\Gamma'_n)$ with weights $0$ and $1$.
\end{lemma}

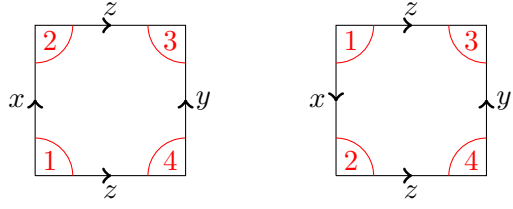
\begin{figure}[ht]
\centering

\begin{tikzpicture}

\draw(0,0)--(2,0)--(2,2)--(0,2)--(0,0);
\draw[->,line width=1.2pt] (1,0)--(1.01,0); \draw[->,line width=1.2pt] (2,1)--(2,1.01); \draw[->,line width=1.2pt] (1,2)--(1.01,2); \draw[->,line width=1.2pt] (0,1)--(0,1.01);
\node [below] at (1,0) {$z$}; \node [above] at (1,2) {$z$}; \node [right] at (2,1) {$y$}; \node [left] at (0,1) {$x$};
\draw [red ] (0.5,0)  arc (0:90:0.5); \draw[red] (0,1.5) arc (-90:0:0.5); \draw[red] (1.5,0) arc (180:90:0.5); \draw[red] (2,1.5) arc (-90:-180:0.5);
\node [red] at (0.2,0.2) {$1$}; \node [red] at (0.2,1.8) {$2$}; \node [red] at (1.8,1.8) {$3$}; \node [red] at (1.8,0.2) {$4$};

\draw (4,0)--(6,0)--(6,2)--(4,2)--(4,0);
\draw[->,line width=1.2pt] (5,0)--(5.01,0); \draw[->,line width=1.2pt] (5,2)--(5.01,2); \draw[->,line width=1.2pt] (4,1)--(4,0.99); \draw[->,line width=1.2pt] (6,1)--(6,1.01); 
\node [below] at (5,0) {$z$}; \node [above] at (5,2) {$z$}; \node [right] at (6,1) {$y$}; \node [left] at (4,1) {$x$};
\draw [red ] (4.5,0)  arc (0:90:0.5); \draw[red] (4,1.5) arc (-90:0:0.5); \draw[red] (5.5,0) arc (180:90:0.5); \draw[red] (6,1.5) arc (-90:-180:0.5);
\node [red] at (4.2,0.2) {$2$}; \node [red] at (4.2,1.8) {$1$}; \node [red] at (5.8,1.8) {$3$}; \node [red] at (5.8,0.2) {$4$};
\end{tikzpicture}
\caption{\label{fig:reorient} $lk(K(\Gamma))=lk(K(\Gamma_x))$: The corners of the original 2-cell also appear in the 2-cell with the edge $x$ reversed, only the order in which the corners appear changes. }
\end{figure}

Before we give a proof we introduce some useful notation. A {\em signed LOT $\Gamma$ } is a labeled oriented tree where we allow vertices to carry signs. An oriented edge from $x^{\epsilon_1}$ to $y^{\epsilon_2}$, $\epsilon_j=\pm 1$, labeled by $z$ gives a 2-cell in $K(\Gamma)$ with attaching map $x^{\epsilon_1}z(zy^{\epsilon_2})^{-1}$. Given a labeled oriented tree $\Gamma$ and a subset $X$ of the vertices of $\Gamma$, we define $\Gamma_{X}$ to be the signed LOT obtained from $\Gamma$ by replacing each vertex $x\in X$ with $x^{-1}$. It is important to note that edge labels and edge orientations in $\Gamma$ and $\Gamma_X$ are the same. One of the key observations is that the link does not change under this vertex sign change: $lk(K(\Gamma))=lk(K(\Gamma_X))$. See Figure \ref{fig:reorient}. In particular if $K(\Gamma)$ satisfies the weight test, then so does $K(\Gamma_X)$.\\

\noindent Proof of Lemma \ref{lem:weights}. Assume $\Gamma'$ is obtained from $\Gamma$ by reversing a single edge labeled $x$. Suppose first that $x$ is not contained in any of the $\Gamma_i$. Then 
\[ lk(K(\Gamma), K(\Gamma_1\cup\ldots\cup \Gamma_n))=lk(K(\Gamma_x), K(\Gamma_1\cup\ldots\cup \Gamma_n)).\]
Since $(K(\Gamma), K(\Gamma_1\cup\ldots\cup \Gamma_n))$ satisfies the weight test with weights $0,1$, so does\\
$(K(\Gamma_x), K(\Gamma_1\cup\ldots\cup \Gamma_n))$. Let $\phi_x\colon K(\Gamma_x)\to K(\Gamma')$ be the homeomorphism that changes the orientation of all $x$-edges in the attaching maps of 2-cells. It induces a homeomorphism of the corresponding expansions and preserves the weight test. Thus $K(\Gamma')$  satisfies the weight test relative $K(\Gamma_1\cup\ldots\cup\Gamma_n)$ with weights $0,1$. 

Next assume that $x$ is contained in one of the $\Gamma_i$, say $\Gamma_1$. If we proceed as above we run into a technical difficulty: the attaching maps of the 2-cells in the subcomplex $K(\Gamma_{1_x}\cup \Gamma_2\cup\ldots\cup \Gamma_n)$ do not all have exponent sum zero, so we are not in the setting of the weight test as given in Definition \ref{def:wtest} anymore. Here is how we fix this. Let $X$ be the set of vertices of $\Gamma_1$. Note that the attaching maps in $K(\Gamma_{1_X})$ do have exponent sum zero. Now we argue exactly as above with $X$ in place of $x$. The homeomorphism $\phi_X\colon K(\Gamma_X)\to K(\Gamma')$ now changes the orientation of all $x$-edges, $x\in X$, in the attaching maps of 2-cells. \qed\\

\begin{prop}\label{prop:case1} Let $\Gamma $ be a reduced injective LOT which is not prime. Assume there is a complete set of sub-LOTs $\Gamma_1,\ldots ,\Gamma_n$. Then $K(\Gamma)$ is VA relative to $K(\Gamma_1)\cup\ldots\cup K(\Gamma_n)$.
\end{prop}

\noindent Proof: By Lemma \ref{lem:reor} there exists a reorientation $\Gamma'$ so that $K(\Gamma')$  satisfies the weight test relative $K(\Gamma'_1\cup\ldots\cup \Gamma'_n)$ with weights $0$ and $1$. By Lemma \ref{lem:weights} $K(\Gamma)$ itself satisfies the weight test relative  $K(\Gamma_1\cup\ldots\cup \Gamma_n)$ with weights $0$ and $1$. It follows from Theorem \ref{thm:WTD} that $K(\Gamma)$ is VA relative to $K(\Gamma_1)\cup\ldots\cup K(\Gamma_n)$.\qed \ \\

\begin{figure}[htbp]
   \centering
   \begin{tikzpicture}

\draw(0,0)--(6,0);
\node [below] at (0,0) {$x_1$}; \fill (0,0) circle (2pt);
\node [below] at (1,0) {$x_2$}; \fill (1,0) circle (2pt);
\node [below] at (2,0) {$x_3$}; \fill (2,0) circle (2pt);
\node [below] at (3,0) {$z$}; \fill (3,0) circle (2pt);
\node [below] at (4,0) {$y_3$}; \fill (4,0) circle (2pt);
\node [below] at (5,0) {$y_2$}; \fill (5,0) circle (2pt);
\node [below] at (6,0) {$y_1$}; \fill (6,0) circle (2pt);
\node [above] at (0.5,0) {$x_3$};
\node [above] at (1.5,0) {$x_1$};
\node [above] at (2.5,0) {$x_2$};
\node [above] at (3.5,0) {$y_2$};
\node [above] at (4.5,0) {$y_1$};
\node [above] at (5.5,0) {$y_3$};

   \end{tikzpicture}
\caption{The figure shows a reduced injective LOT that does not contain a complete set of sub-LOTs (for any orientation of its edges). Note that it freely decomposes.}
   \label{fig:nocomp}
\end{figure}
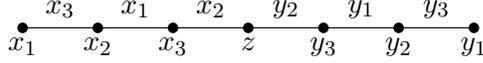

This proposition would quickly lead to a proof of our main Theorem \ref{thm:general} if all non-prime reduced injective LOTs would have a complete set of sub-LOTs. This, however, is not true. See Figure \ref{fig:nocomp}. We thank the referee for catching this problem.\\

We say a LOT $\Gamma$ {\em freely decomposes} if $\Gamma=\Gamma_L\cup \Gamma_R$, where the $\Gamma_L$ and $\Gamma_R$ are proper sub-LOTs of $\Gamma$, intersecting in a single vertex.

\vfill\eject
\begin{prop}\label{prop:comp} Suppose $\Gamma$ is a reduced and injective LOT that is not prime. Then either
\begin{enumerate}
\item $\Gamma$ contains a complete set of sub-LOTs, or 
\item $\Gamma$ freely decomposes.
\end{enumerate}
\end{prop}

\noindent Proof. We make two observations:\\

\noindent (1) {\em Let $\Gamma$ be an injective compressed LOT and let $\Gamma_1$ be a maximal proper sub-LOT. If the quotient $\bar \Gamma_1$ is not compressed, then $\Gamma=\Gamma_1\cup e$, where $e$ is an edge attached to $\Gamma_1$, whose vertex $y$ not in $\Gamma_1$ does not occur as an edge label in $\Gamma_1$. In particular, $\Gamma$ is not boundary reduced.}\\

There exists an edge $\bar e$ in $\bar \Gamma_1$ with vertices $\bar x$ and $\bar y$ and edge label $\bar x$. Since $\Gamma$ is compressed it follows that $\bar x$ is the vertex the sub-LOT $\Gamma_1$ is collapsed to, and it follows that we have a sub-LOT $\Gamma_1\cup e$ in $\Gamma$, where the edge label of $e$ is a vertex in $\Gamma_1$. So $\Gamma_1\cup e$ is a sub-LOT of $\Gamma$ containing $\Gamma_1$. Maximality of $\Gamma_1$ in $\Gamma$ implies that  $\Gamma=\Gamma_1\cup e$. All vertices in $\Gamma_1$ except $x$ occur in $\Gamma_1$ as edge labels (by the way we collapse), so $y$ is not an edge label of $\Gamma_1$. \\

\noindent (2) {\em Let $\Gamma$ be an injective compressed LOT and suppose 
   $$\begin{matrix} 
     \Gamma & \to &  \bar \Gamma_1& \to &\bar \Gamma_{12}  \\
       \cup      &       &   \cup                       &       &                                                                \\
     \Gamma_1&      & \bar \Gamma_2 &
       \end{matrix}$$ 
is a 2-stage quotient of $\Gamma$ so that $\bar \Gamma_1$ is compressed but $\bar \Gamma_{12}$ is not. Then $\Gamma$ freely decomposes.}\\

We know from Observation (1) that $\bar \Gamma_1=\bar \Gamma_2\cup \bar e$ where $\bar e$ is attached to $\bar \Gamma_2$ at a vertex $\bar x$, the other vertex $\bar y$ does not label an edge in $\bar \Gamma_2$, and $\bar e$ has as edge label a vertex from $\bar \Gamma_2$. Let $\bar z$ be the vertex $\Gamma_1$ is collapsed to. Note that $\bar z$ can not be in $\bar \Gamma_2$ because of maximality, so $\bar z=\bar y$. But then it follows that $\Gamma=\Gamma_2\cup e\cup\Gamma_1$. If we take $\Gamma_L=\Gamma_2\cup e$ and $\Gamma_R=\Gamma_1$ we have a free decomposition of $\Gamma$. \\

Using the two observations we complete the proof of Proposition \ref{prop:comp}. Suppose $\Gamma$ is not prime and $\bar\Gamma_{12...n}$ is an injective prime quotient that is not compressed. Note that $n\ge 2$ because of the first observation (1): if $n=1$ then $\Gamma$ is not boundary reduced, contradicting the assumption that $\Gamma$ is reduced. 
There exists a $j$ so that $\bar \Gamma_{12...j}\to  \bar \Gamma_{12...j+1}\to \bar \Gamma_{12...j+2}$, so that both $\bar \Gamma_{12...j}$ and $\bar \Gamma_{12...j+1}$ are injective and compressed, but $\bar \Gamma_{12...j+2}$ is not. Then by observation (2) $\bar \Gamma_{12...j}$ freely decomposes: $\bar \Gamma_{12...j}=\bar \Gamma_L\cup \bar \Gamma_R$.
We claim $\Gamma =\bar \Gamma_{12...j}$. If not, then let $\bar x$ be the vertex in $\bar \Gamma_{1...j}$ that $\bar \Gamma_j\subseteq \bar \Gamma_{12...(j-1)}$ is collapsed to. But then $\bar x$ is in $\bar \Gamma_L$ or in $\bar \Gamma_R$, contradicting maximality of $\bar \Gamma_j$ in $\bar \Gamma_{12...(j-1)}$.
So  $\Gamma =\bar \Gamma_{12...j}$,  and we have that $\Gamma$ freely decomposes. \qed \\

We need two more lemmas before we can complete the proof of Theorem \ref{thm:general}.

\begin{lemma}\label{lem:bdrred} Suppose $\Gamma_0\subseteq \Gamma$ and $\Gamma_0$ is obtained from $\Gamma $ by a boundary reduction. Then $K(\Gamma)$ is VA if and only if $K(\Gamma_0)$ is VA.
\end{lemma} 

\noindent Proof: If $K(\Gamma)$ is VA, then so is $K(\Gamma_0)$ because the VA property is hereditary. Assume that $K(\Gamma_0)$ is VA. If $\Gamma$ consists of one edge then $K(\Gamma)$ is VA. Let $\Gamma$ contain at least two edges. Then $\Gamma=\Gamma_0\cup e$ where $e$ is attached to $\Gamma_0$ at a vertex $x$, and the other vertex $y$ of $e$ does  not occur as an edge label in $\Gamma$. $K(\Gamma )$ contains a 2-cell with a free boundary edge $y$. So any vertex reduced spherical diagram $f\colon C\to K(\Gamma )$ maps to $K(\Gamma_0)$.\qed\\


\begin{lemma}\label{lem:compVA} If an injective LOT $\Gamma$ freely decomposes $\Gamma=\Gamma_L\cup \Gamma_R$ and both $K(\Gamma_L)$ and $K(\Gamma_R)$ are VA, then $K(\Gamma)$ is VA.
\end{lemma}

\noindent Proof: Consider a vertex reduced spherical diagram $f\colon C\to K(\Gamma)$. $f(C)$ cannot be contained in one of $K(\Gamma_L)$ or $K(\Gamma_R)$ only, because $K(\Gamma_L)$ and $K(\Gamma_R)$ are VA. Let $x$ be the common vertex of $\Gamma_L$ and $\Gamma_R$. $x$ is edge label at most once in $\Gamma$ since $\Gamma$ is injective. So assume $x$ is not edge label in $\Gamma_R$. Let $\gamma \in C$ be a boundary component of a maximal region, which maps to $K(\Gamma_R)$. Each edge of $\gamma$ maps to $x$, and so $\gamma$ maps to a word in $x$ and $x^{-1}$ of exponent sum zero. So there is a vertex $v\in \gamma$ where two edges $e_1,e_2\in \gamma$ end. Now $v$ is a folding vertex with folding pair of 2-cells $d_1$ and $d_2$ which both map to $K(\Gamma_R)$ containing $e_1$ and $e_2$, respectively.\qed \\

{\bf Proof of Theorem \ref{thm:general}}: Let $\Gamma$ be a smallest injective compressed LOT such that $K(\Gamma)$ is not VA. Minimality implies that $\Gamma $ is boundary reduced because of Lemma \ref{lem:bdrred}, so $\Gamma$ is reduced.
$\Gamma $ is not prime because of Theorem \ref{thm:prime}. $\Gamma$ does not freely decompose by Lemma \ref{lem:compVA}. So by Proposition \ref{prop:comp} $\Gamma$ contains a complete set of sub-LOTs $\Gamma_1,\ldots ,\Gamma_n$. Proposition \ref{prop:case1} implies that $K(\Gamma)$ is VA relative to $K(\Gamma_1)\cup\ldots\cup K(\Gamma_n)$. Since $K(\Gamma)$ is the smallest non-VA injective compressed LOT we have that all the $K(\Gamma_i)$ are VA. But then Theorem \ref{sVAsub} implies that $K(\Gamma)$ is VA in contradiction to our assumption. \qed


\end{document}